\documentclass[12pt]{article}
\usepackage{graphicx,amsmath,amssymb}

\def \R {{\mathbb R}}
\def \Z {{\mathbb Z}}
\def \N {{\mathbb N}}
\def \T {{\mathbb T}}
\def \Z {{\mathbb Z}}
\def \x {{\bf{x}}}
\def \y {{\bf{y}}}
\def \a {{\bf{a}}}

\def \n {{\bf{n}}}
\def \V {{\bf{V}}}
\def \D {\nabla}
\def \S {{\mathcal S}}
\newcommand{\disp}{\displaystyle}

\newcommand{\be}{\begin{equation}}
\newcommand{\ee}{\end{equation}}

\begin{document}

\normalsize 
\begin{center}
\Large \bf Synchronized Front Propagation and Delayed Flame Quenching in Strain G-equation and Time-Periodic Cellular Flows
\end{center}
\vspace*{0.5cm}

\begin{center}\bf 
Yu-Yu Liu
\footnote{Department of Mathematics, 
National Cheng Kung University, 
Tainan 70101, Taiwan. 
E-mail: yuyul@ncku.edu.tw},
Jack Xin 
\footnote{Department of Mathematics, 
University of California, Irvine, CA 92697, USA.
Email: jxin@math.uci.edu}
\end{center}
\vspace*{0.5cm}

\noindent {\bf Abstract.}
G-equations are level-set type Hamilton-Jacobi partial differential equations modeling propagation of flame front along a flow velocity and a laminar velocity. In consideration of flame stretching, strain rate may be added into the laminar speed. We perform finite difference computation of G-equations with the discretized strain term  being monotone with respect to one-sided spatial derivatives. Let the flow velocity be the time-periodic cellular flow (modeling Rayleigh-Bénard advection), we compute the turbulent flame speeds as the asymptotic propagation speeds from a planar initial flame front. In strain G-equation model, front propagation is enhanced by the cellular flow, and flame quenching occurs if the flow intensity is large enough. In contrast to the results in steady cellular flow, front propagation in time periodic cellular flow may be locked into certain spatial-temporal periodicity pattern, and turbulent flame speed becomes a piecewise constant function of flow intensity. Also the disturbed flame front does not cease propagating until much larger flow intensity.
\bigskip

\bigskip


\newpage

\section{Introduction}
\setcounter{equation}{0}

Front propagation in turbulent combustion is a complex multiscale dynamical process. To analyze and measure the turbulent burning velocity is of great importance in both combustion theory and experiment. Important issues include: (i) front speed enhancement by flow velocity, (ii) bending of front speed growth in large flow velocity and (iii) flame quenching due to flame stretching \cite{B92,MK99,R95,VCK03,X00,X09}.  

In this paper we consider the {\it inviscid G-equation} \cite{P00,W85}
\be\label{Gi}
{\partial G\over\partial t}+\V(\x,t)\cdot\D G+s_L|\D G|=0,
\ee
and the {\it strain G-equation}
\be\label{Gs}
{\partial G\over\partial t}+\V(\x,t)\cdot \D G+s_L|\D G|
+d_M{\D G\cdot D\V\cdot \D G\over|\D G|}=0
\ee
in two-dimensional space ($\x=\left\langle x,y\right\rangle\in \R^2, t>0$). In the corrugated flamelet regime of premixed turbulent combustion, the flame front is considered as the interface $\{G(\x,t)=0\}$ between the burnt region $\{G<0\}$ and the unburnt region $\{G>0\}$. The motion law of flame front
\be\label{law}
{d\x\over dt}=\V(\x,t)+s_L\n
\ee
consists of a prescribed flow velocity $\V(\x,t)$ and a  laminar velocity normal to the level set  $\n=\D G/|\D G|$ with {\it laminar flame speed} $s_L>0$. This motion law gives the inviscid G-equation (\ref{Gi}). To cooperate the flame stretching effect by flow velocity, a correction term may be added into the laminar speed  \cite{LXY13a,P00}:
$$
\hat{s}_L=s_L-d_M\mathcal{S},
$$
where $\mathcal{S}=-\n^t\!\cdot\!D\V\!\cdot\!\n$ is the strain rate with {\em Markstein diffusivity} $d_M>0$. This modified motion law gives the strain G-equation (\ref{Gs}).

\begin{figure}\center
\includegraphics[width=0.4\textwidth]{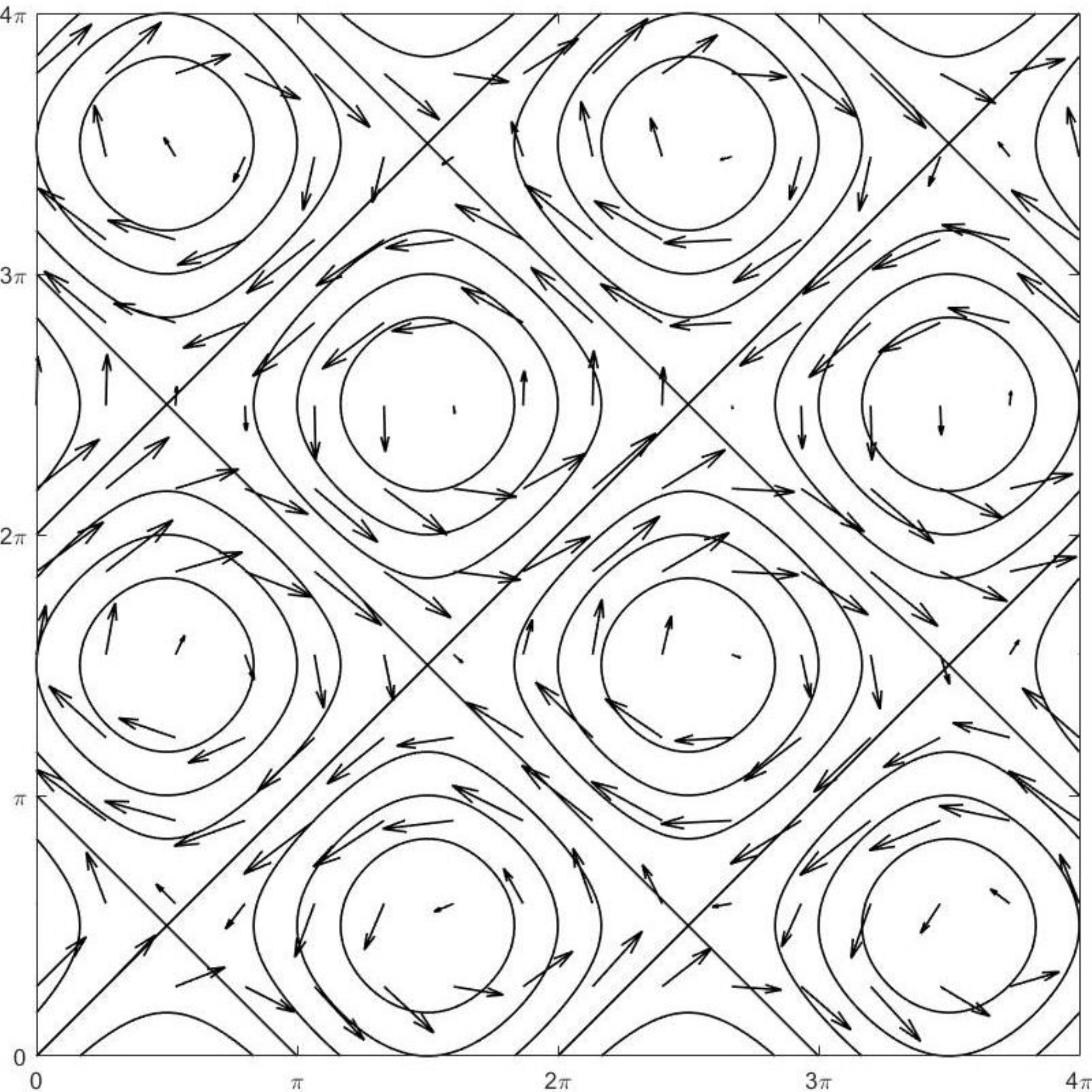}
\includegraphics[width=0.4\textwidth]{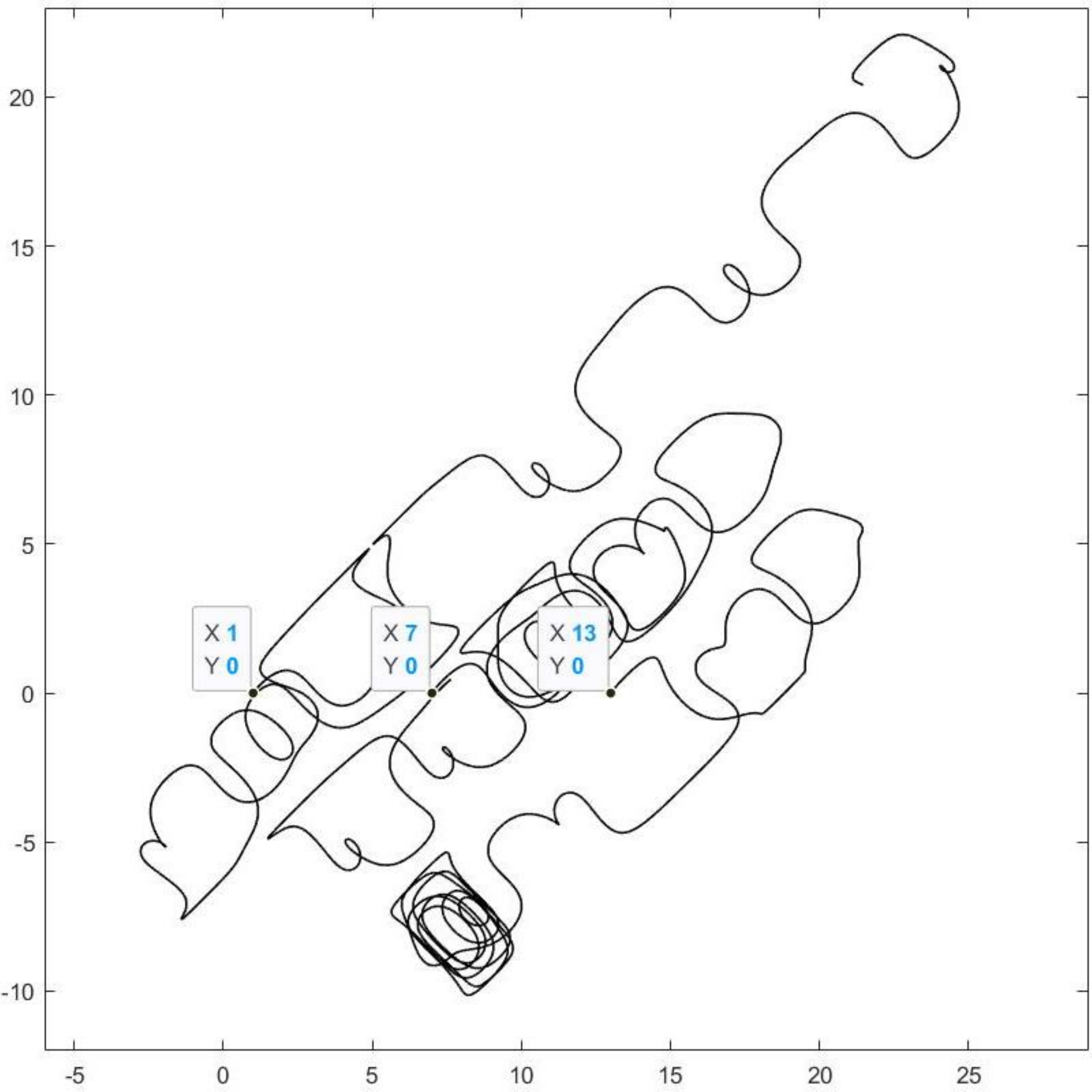}
\caption{Left panel: steady cellular flow (\ref{Vs}). Right panel: trajectories of unsteady cellular flow (\ref{Vu}) with $A=1$, $\omega=1$ and time up to 100.}
\label{Fig_Cell}
\end{figure}

For the flow velocity, we consider the steady cellular flow
\be\label{Vs}
\V(x,y)=A\cdot\left\langle\cos(y),\cos(x)\right\rangle
\ee
and the Rayleigh-Bénard advection \cite{CW91}
\be\label{Vu}
\V(x,y,t)=A\cdot\left\langle
\cos(y)+\sin(y)\cos(\omega t),
\cos(x)+\sin(x)\cos(\omega t)
\right\rangle,\ee
where $A$ is the flow intensity. The Rayleigh-Bénard advection is an unsteady cellular flow periodic in time upon rewriting in the form:
$$\V(x,y,t)=A\cdot\sec(\theta_{\omega}(t))\cdot\left\langle
\cos(y+\theta_{\omega}(t)),\cos(x+\theta_{\omega}(t))
\right\rangle$$
with $\theta_{\omega}(t)=\tan^{-1}(\cos(\omega t))$. Also the Rayleigh-Bénard advection is known for chaotic streamlines and diffusion-like transport in diagonal direction \cite{BCVV95,ZCX15}. See figure \ref{Fig_Cell}. 

If the initial flame front is planar and the flow velocity is at rest ($A=0$), then the flame front propagates at speed $s_L$. If the flow velocity is in motion ($A>0$), then the flame front is corrugated in time and eventually propagates at an asymptotic speed $s_T$ called the {\it turbulent flame speed}. Our goal is to study the growth of turbulent flame speed with respect to the increase of flow intensity. Specifically, we would like to see the qualitative difference of function $s_T(A)$ between the G-equations (\ref{Gi}) and (\ref{Gs}) as well as the cellular flows (\ref{Vs}) and (\ref{Vu}).

The inviscid G-equation (\ref{Gi}) and cellular flow (\ref{Vs}) have been studied in many contexts. In \cite{CNS11,XY10}, formulation of turbulent flame speed is rigorously justified by periodic homogenization theory. In \cite{ACVV02,CTVV03}, motion of level set is simulated by the motion law (\ref{law}) with grid points labeled as burnt or unburnt particles. In \cite{O00,XY13}, estimates of front speed enhancement is obtained using the optimal control  representation of solutions. The problem becomes much more challenging if the strain term is added. In \cite{LXY13a,LXY13b}, computational study of strain G-equation is given along with the curvature term. In \cite{XY14}, flame quenching in steady cellular flow (\ref{Vs}) is rigorously justified using the differential game representation of solutions. 

In \cite{LXY13a} we consider the full G-equation model with motion law $$\hat{s}_L=s_L-d_M(\S+s_L\kappa),$$ where $\kappa=\mbox{div}(\n)$ is the curvature of level set added as the flame stretching effect by the laminar velocity. The curvature term is a nonlinear diffusion that brings certain smoothness to the solution. In the framework of monotone discretization of finite difference computation in Hamilton-Jacobi equations \cite{CL84}, we evaluate the strain rate $\S$ and the curvature $\kappa$ by central differencing in order to apply the Godunov scheme on $\hat{s}_L|\nabla G|$. If the curvature term is removed as in present paper, solutions of strain G-equation (\ref{Gs}) is literally not differentiable. Therefore we shall construct a numerical Hamiltonian of the strain term in (\ref{Gs}) that is genuinely monotone with respect to all one-sided spatial derivatives of the solutions.

In \cite{LXY13b} we consider the inviscid G-equation and the cellular flows that are unsteady in $x$-direction (see also \cite{CTVV03}), and we evaluate turbulent flame speeds as a function of the temporal frequency: $s_T(\omega)$. It is observed that front propagation may be synchronized with spatial and temporal periodicity of the cellular flow. Then turbulent flame speed is a piecewise linear function of frequency with rational slopes: $s_T(\omega)=r\omega$, $\omega\in I_r$ for some $r\in\mathbb{Q}$ and intervals $I_r\subset\R$. In this paper we will consider $s_T(A)$ with the strain term being added.

The rest of the paper is organized as follows. In section 2, we construct the numerical discretization of G-equation models. In section 3, we present the numerical results of  turbulent flame speeds. In section 4,  we conclude the paper with future work and acknowledgments.

\section{Numerical Hamiltonian of G-equations}

The general form of Hamilton-Jacobi (HJ) equations are
\be\label{HJ}
{\partial G\over\partial t}+H(D_xG,D_yG)=0
\ee
where $G(x,y,t):\R^2\times\R\to\R$ is the solution and $H(p,q):\R^2\to\R$ is the Hamiltonian. The solutions are defined in viscosity sense and may not be differentiable. Let the uniform discretization of the solutions be $G^n_{i,j}=G(i\Delta x,j\Delta y,n\Delta t)$, then the finite difference and forward Euler discretization of (\ref{HJ}) is
$$
{G^{n+1}_{i,j}-G^n_{i,j}\over\Delta t}+\hat{H}(
D_x^-G^n_{i,j},D_x^+G^n_{i,j},
D_y^-G^n_{i,j},D_y^+G^n_{i,j})=0,
$$
where $D_x^-G^n_{i,j}$, $D_y^-G^n_{i,j}$, $D_y^-G^n_{i,j}$ and $D_y^+G^n_{i,j}$ are one-sided approximations of the spatial derivatives and $\hat{H}(p^-,p^+,q^-,q^+)$ is the numerical Hamiltonian of $H(p,q)$. To obtain the numerical stability, $\hat{H}$ is chosen to be consistent ($\hat{H}(p,p,q,q)=H(p,q)$) and monotone (symbolically $\hat{H}(\uparrow,\downarrow,\uparrow,\downarrow)$). A popular choice is the Lax-Friedrichs scheme:
$$\textstyle
\hat{H}(p^-,p^+,q^-,q^+)=H({p^-+p^+\over2},{q^-+q^+\over2})-\|{\partial H\over\partial p}\|_{\infty}({p^+-p^-\over2})-\|{\partial H\over\partial q}\|_{\infty}({q^+-q^-\over2}).
$$
But it is desirable to reduce the artificial diffusion $(p^+-p^-)/2$, $(q^+-q^-)/2$ whenever possible.

To improve the accuracy of the solutions, the spatial derivatives are evaluated by high order WENO (weighted essentially non-oscillatory) scheme, and the time steps are iterated by high order TVD (total variation diminishing) Runge-Kutta (RK) scheme. A popular choice is fifth order scheme in space (WENO5) paired with third order scheme in time (TVD-RK3). Time step size $\Delta t$ is determined by the CFL (Courant–Friedrichs–Lewy) condition. Here we present the construction of numerical Hamiltonian for inviscid and strain G-equations (\ref{Gi})(\ref{Gs}) and refer \cite{JP00,OF02,S07} for implementation of WENO and TVD-RK schemes. 

For inviscid G-equation (\ref{Gi}), write $\V=\left<u,v\right>$ and the Hamiltonian is
$$H_{\mathrm{inv}}(p,q)=up+vq+s_L\sqrt{p^2+q^2}.$$
The corresponding numerical Hamiltonian is
$$\hat{H}_{\mathrm{inv}}(p^-,p^+,q^-,q^+)=up_{\mathrm{vel}}+vq_{\mathrm{vel}}
+s_L\sqrt{p_{\mathrm{lem}}^2+q_{\mathrm{lem}}^2},$$
where the spatial derivatives in the velocity term are given by the upwind scheme:
$$p_{\mathrm{vel}}=\left\{\!\begin{array}{ll}
p^- & ,\mathrm{if}\ u>0\\
p^+ & ,\mathrm{if}\ u<0
\end{array}\right.,
q_{\mathrm{vel}}=\left\{\!\begin{array}{ll}
q^- & ,\mathrm{if}\ v>0\\
q^+ & ,\mathrm{if}\ v<0
\end{array}\right.,$$
and the spatial derivatives in the laminar term are given by the Godunov scheme:
$$p_{\mathrm{lem}}^2=\max(\max(p^-,0)^2,\min(p^+,0)^2),$$
$$q_{\mathrm{lem}}^2=\max(\max(q^-,0)^2,\min(q^+,0)^2).$$

For strain G-equation (\ref{Gs}), it suffices to consider the strain term with Hamiltonian
\be\label{Hs}
H_{\mathrm{str}}(p,q)=a{p^2\over\sqrt{p^2+q^2}}
+b{q^2\over\sqrt{p^2+q^2}}+c{pq\over\sqrt{p^2+q^2}},
\ee
where $a=d_M(\partial u/\partial x)$, $b=d_M(\partial v/\partial y)$ and $c=d_M(\partial u/\partial y+\partial v/\partial x))$. 

The first term of (\ref{Hs})
$$
H_{\mathrm{s1}}(p,q)=a{p^2\over\sqrt{p^2+q^2}}
$$
is monotone increasing with respect to $p^2$ and monotone decreasing with respect to $q^2$ if $a>0$ (opposite monotonicity if $a<0$). Therefore its numerical Hamiltonian is given by the Osher-Sethian scheme:
$$
\hat{H}_{\mathrm{s1}}(p^-,p^+,q^-,q^+)
=a{p_{\mathrm{s1}}^2\over\sqrt{p_{\mathrm{s1}}^2+q_{\mathrm{s1}}^2}},
$$
$$
p_{\mathrm{s1}}^2=\left\{\!\begin{array}{ll}
\min(p^+,0)^2+\max(p^-,0)^2 & ,\mathrm{if}\ a>0\\
\min(p^-,0)^2+\max(p^+,0)^2 & ,\mathrm{if}\ a<0
\end{array}\right.,
$$
$$
q_{\mathrm{s1}}^2=\left\{\!\begin{array}{ll}
\min(q^-,0)^2+\max(q^+,0)^2 & ,\mathrm{if}\ a>0\\
\min(q^+,0)^2+\max(q^-,0)^2 & ,\mathrm{if}\ a<0
\end{array}\right..
$$

The second term of (\ref{Hs})
$$
H_{\mathrm{s2}}(p,q)=b{q^2\over\sqrt{p^2+q^2}}
$$
is monotone decreasing with respect to $p^2$ and monotone increasing with respect to $q^2$ if $b>0$ (opposite monotonicity if $b<0$). Therefore its numerical Hamiltonian is given by the Osher-Sethian scheme:
$$
\hat{H}_{\mathrm{s2}}(p^-,p^+,q^-,q^+)
=b{p_{\mathrm{s2}}^2\over\sqrt{p_{\mathrm{s2}}^2+q_{\mathrm{s2}}^2}},
$$
$$
p_{\mathrm{s2}}^2=\left\{\!\begin{array}{ll}
\min(p^-,0)^2+\max(p^+,0)^2 & ,\mathrm{if}\ b>0\\
\min(p^+,0)^2+\max(p^-,0)^2 & ,\mathrm{if}\ b<0
\end{array}\right.,
$$
$$
q_{\mathrm{s2}}^2=\left\{\!\begin{array}{ll}
\min(q^+,0)^2+\max(q^-,0)^2 & ,\mathrm{if}\ b>0\\
\min(q^-,0)^2+\max(q^+,0)^2 & ,\mathrm{if}\ b<0
\end{array}\right..
$$

Finally consider the third term of (\ref{Hs})
$$
H_{\mathrm{s3}}(p,q)=c{pq\over\sqrt{p^2+q^2}}.
$$
Observe that
$${\partial H_{\mathrm{s3}}\over\partial p}
={cq^3\over(p^2+q^2)^{3\over2}}, 
{\partial H_{\mathrm{s3}}\over\partial q}
={cp^3\over(p^2+q^2)^{3\over2}}.$$
Then upwind direction of $p$ is determined if $q^-$, $q^+$ have same sign, and upwind direction of $q$ is determined if $p^-$, $p^+$ have same sign. Also
$\partial H_{\mathrm{s3}}/\partial p>0$ if $cq>0$,
$\partial H_{\mathrm{s3}}/\partial p<0$ if $cq<0$,
$\partial H_{\mathrm{s3}}/\partial q>0$ if $cp>0$ and
$\partial H_{\mathrm{s3}}/\partial q<0$ if $cp<0$. Otherwise Lax-Friedrichs scheme is applied with
$|\partial H_{\mathrm{s3}}/\partial p|\leq |c|$ and 
$|\partial H_{\mathrm{s3}}/\partial q|\leq |c|$.
Therefore its numerical Hamiltonian is given by the Roe scheme:
$$
\hat{H}_{\mathrm{s3}}(p^-,p^+,q^-,q^+)=
c{p_{\mathrm{s3}}q_{\mathrm{s3}}\over\sqrt{p_{\mathrm{s3}}^2+q_{\mathrm{s3}}^2}}
-\bar{c}_p{(p^+-p^-)\over 2}-\bar{c}_q{(q^+-q^-)\over 2},
$$
$$p_{\mathrm{s3}}=\left\{\!\begin{array}{ll}
p^-&,\mbox{if $q^-q^+>0$ and $cq^{\pm}>0$}\\
p^+&,\mbox{if $q^-q^+>0$ and $cq^{\pm}<0$}\\
\disp{p^++p^-\over2}& ,\mbox{if $q^-q^+<0$}
\end{array}\right.,$$
$$
q_{\mathrm{s3}}=\left\{\!\begin{array}{ll}
q^-&,\mbox{if $p^-p^+>0$ and $cp^{\pm}>0$}\\
q^+&,\mbox{if $p^-p^+>0$ and $cp^{\pm}<0$}\\
\disp{q^++q^-\over2}& ,\mbox{if $p^-p^+<0$}
\end{array}\right.,$$
$$
\bar{c}_p=\left\{\!\begin{array}{ll}
0& ,\mbox{if $q^-q^+>0$}\\
|c|& ,\mbox{if $q^-q^+<0$}
\end{array}\right.,
\bar{c}_q=\left\{\!\begin{array}{ll}
0& ,\mbox{if $p^+p^->0$}\\
|c|& ,\mbox{if $p^+p^-<0$}
\end{array}\right..$$

Overall the artificial diffusion is added only when the one-sided derivatives have opposite signs.

\section{Numerical Results}

Let $G(\x,0)=x$ for $\x\in\R^2$, then the flame front is initially $\{x=0\}$ and starts propagating in $x$-direction. Note that the cellular flows (\ref{Vs})(\ref{Vu}) are spatially periodic on $(2\pi\T)^2$, we may write $G(\x,t)=x+u(\x,t)$ with $u(\x,t)$ spatially periodic for all $t>0$. Therefore we can solve the initial-boundary value problem of inviscid G-equation (\ref{Gi}) in finite spatial domain:
\be\label{IVP}
\left\{\begin{array}{ll}
{\partial G\over\partial t}+\V(\x,t)\cdot \D G+s_L|\D G|=0&, \x\in[0,2\pi]^2, t>0\\
G(\x,0)=x&, \x\in [0,2\pi]^2\\
G(x,2\pi,t)=G(x,0,t) &, x\in [0,2\pi], t>0\\
G(2\pi,y,t)=G(0,y,t)+2\pi &, y\in [0,2\pi], t>0
\end{array}\right..
\ee
Initial-boundary conditions for strain G-equation (\ref{Gs}) are exactly the same. 

\begin{figure}\center
\includegraphics[width=0.8\textwidth]{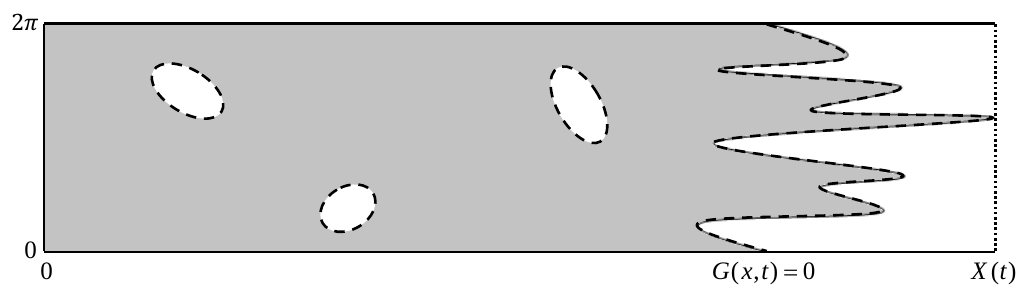}
\caption{Propagation distance $X(t)$.}
\label{Fig_Xt}
\end{figure}

Numerical computation of (\ref{IVP}) is carried out on a $256\times256$ uniform mesh of spatial domain $[0,2\pi]^2$.
Then we can obtain the solution on stripe domain $\R\!\times\![0,2\pi]$ with
$$
G(x+2k\pi,y,t)=G(x,y,t)+2k\pi, k\in\Z
$$
so that we can visualize the level set $\{G(\x,t)=0\}$. Denote the propagation distance of the flame front in $x$-direction in time as follows:
\be\label{Xt}
X(t)=\sup\{x\in\R\,|\,G(\x,t)<0\}.
\ee
See figure \ref{Fig_Xt}. We see $X(0)=0$ and $X'(t)$ is the instantaneous propagation speed. Therefore the turbulent flame speed is defined as the asymptotic propagation speed in large time as follows:
$$
s_T=\lim_{t\to\infty}{X(t)\over t}.
$$

\begin{figure}\center
\includegraphics[width=0.9\textwidth]{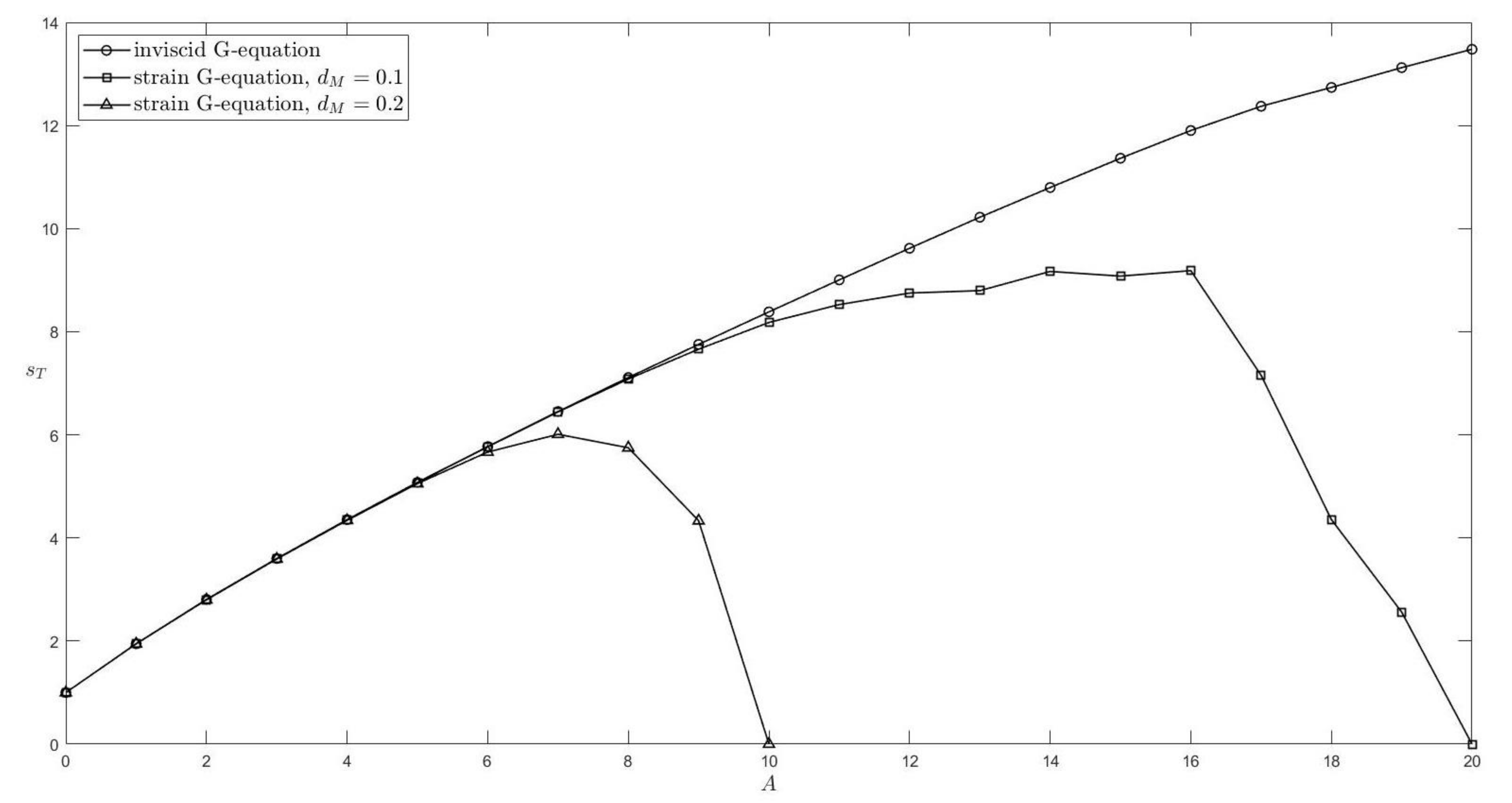}
\caption{Turbulent flame speed $s_T(A)$ for inviscid G-equation (\ref{Gi}) and strain G-equation (\ref{Gs}) with $d_M=0.1, 0.2$ and steady cellular flow (\ref{Vs}).}
\label{Fig_StVsA}
\end{figure}

\begin{figure}\center
\includegraphics[width=0.9\textwidth]{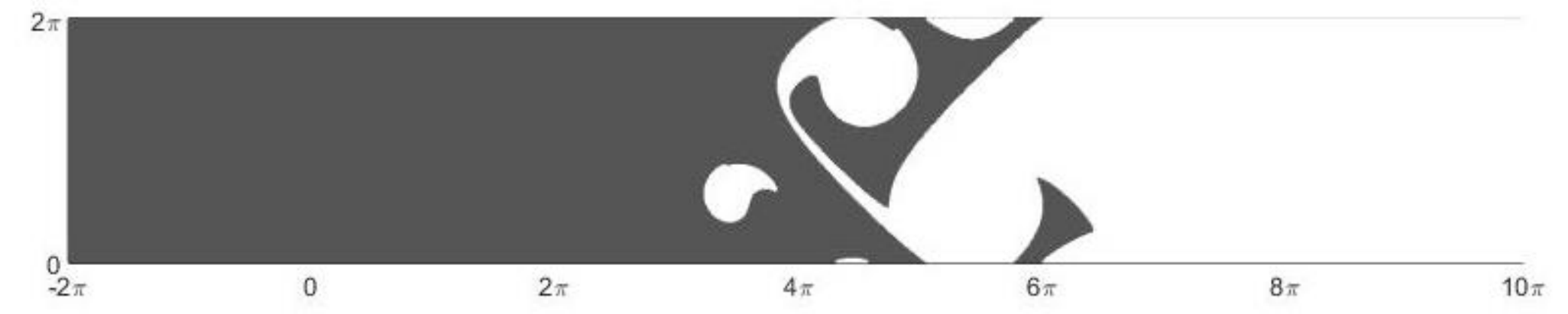}
\includegraphics[width=0.9\textwidth]{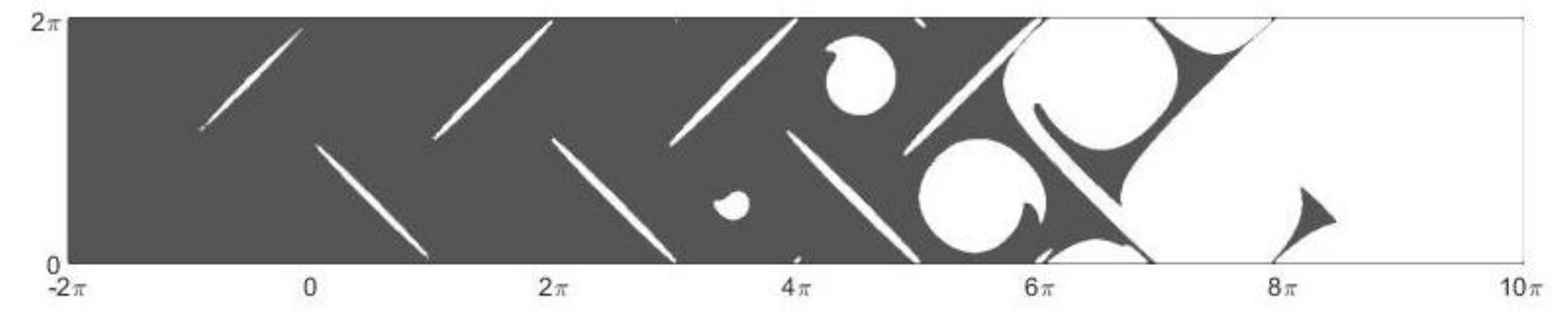}
\includegraphics[width=0.9\textwidth]{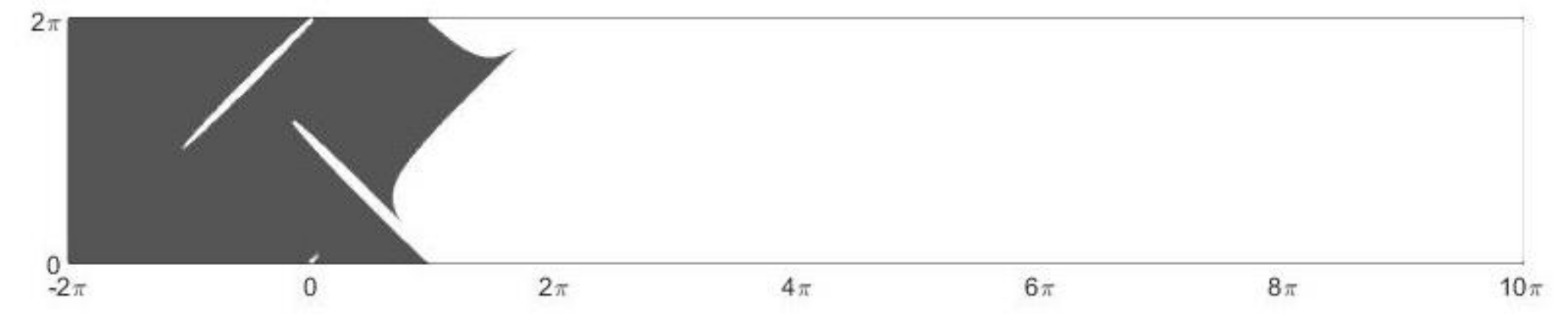}
\caption{Front propagation in strain G-equation (\ref{Gs}) with $d_M=0.2$ and steady cellular flow (\ref{Vs}) at time $t=4$. Upper panel: $A=5$ (complete combustion). Middle panel: $A=9$ (incomplete combustion). Lower panel: $A=12$ (flame quenching).}
\label{Fig_GVsA}
\end{figure}

Figure \ref{Fig_StVsA} shows the plots of $s_T(A)$ of G-equations (\ref{Gi})(\ref{Gs}) with steady cellular flow (\ref{Vs}). For inviscid G-equation, turbulent flame speed is enhanced by the cellular flow with growth rate $s_T=O(A/\log A)$, $A\gg1$. The sublinear growth is due to slowdown of front propagation near the hyperbolic equilibria of the cellular flow. For strain G-equation, turbulent flame speed starts to decrease and soon drops to zero for larger flow intensity. Figure \ref{Fig_GVsA} further presents three stages of front propagation being affected by the strain term as the flow intensity increases. When $A$ is relatively small, the strain rate $\S$ is so small that the laminar speed $\hat{s}_L$ remains strictly positive. Therefore the flame front propagates without unburned region being left behind (complete combustion). When $A$ is moderately larger, the laminar speed $\hat{s}_L$ decreases as the strain rate $\S$ increases near the hyperbolic equilibria. The flame front still manages to propagate forward, but there exists stagnated unburnt regions being left behind (incomplete combustion). When $A$ exceeds a certain value, the strain rate $\S$ is large enough to negate both the flow velocity and the laminar velocity. Therefore the flame front ceases to propagate forward (flame quenching).

\begin{figure}\center
\includegraphics[width=0.8\textwidth]{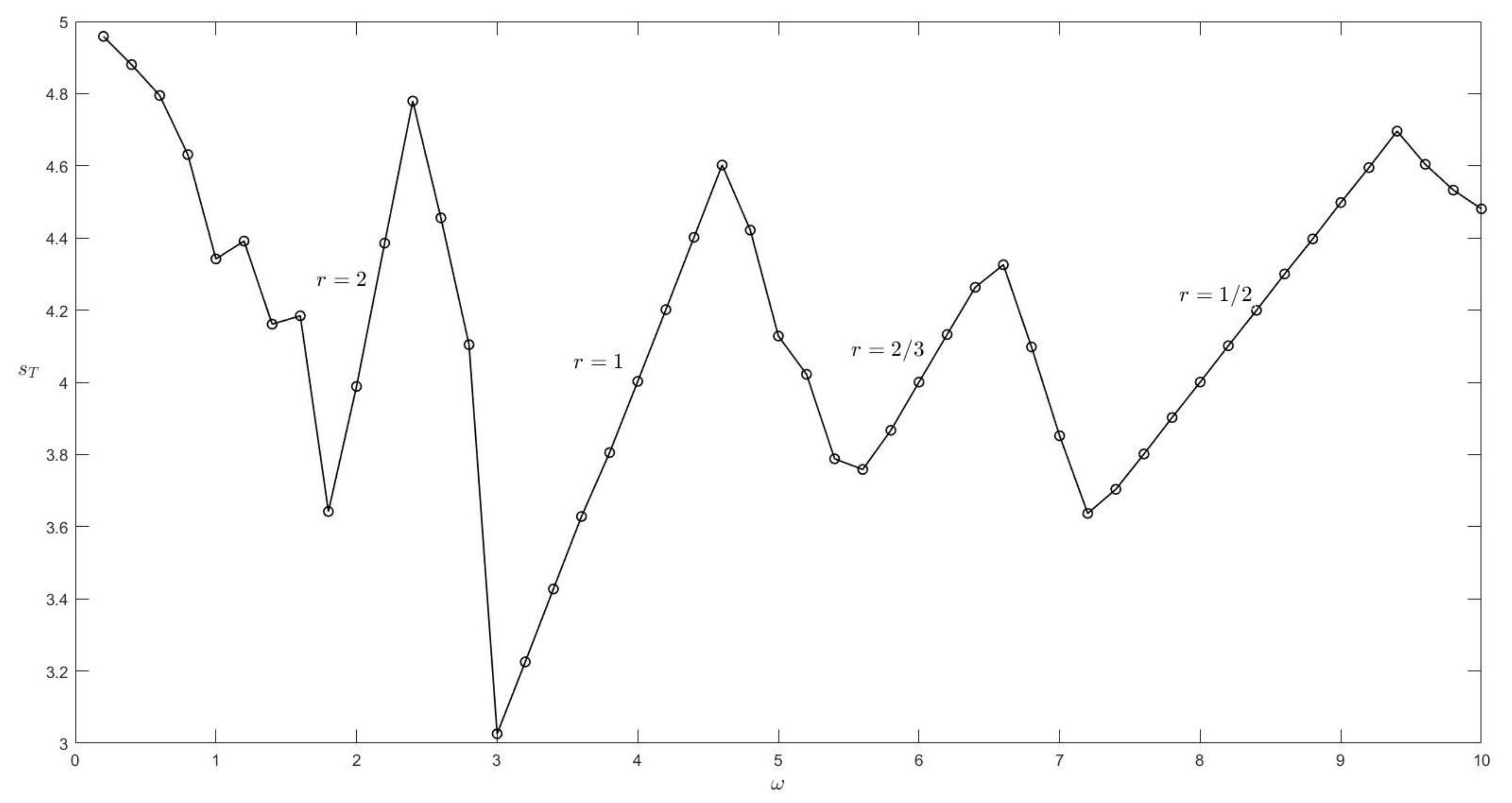}
\caption{Turbulent flame speed $s_T(\omega)$ for inviscid G-equation (\ref{Gi}) and unsteady cellular flow (\ref{Vu}) with $A=4$.}
\label{Fig_StGiVuW}
\end{figure}

Figure \ref{Fig_StGiVuW} shows the plot of $s_T(\omega)$ for inviscid G-equation (\ref{Gi}) and unsteady cellular flow (\ref{Vu}) with $A=4$. It happens that propagation of flame front is eventually synchronized with the spatial and temporal periodicity of the cellular flow. Specifically, the propagation distance is a multiple of the spatial period $\Delta x=2\pi\cdot N$, and the propagation time is a multiple of temporal period $\Delta t=2\pi/\omega\cdot M$. Also the synchronization pattern $N,M\in\N$ is robust with respect to small variation of $\omega$ called {\em frequency locking}. Therefore the turbulent flame speed  $s_T(\omega)=\Delta x/\Delta t=r\omega$ is a piecewise linear function with rational slope $r=N/M$. 

\begin{figure}\center
\includegraphics[width=0.9\textwidth]{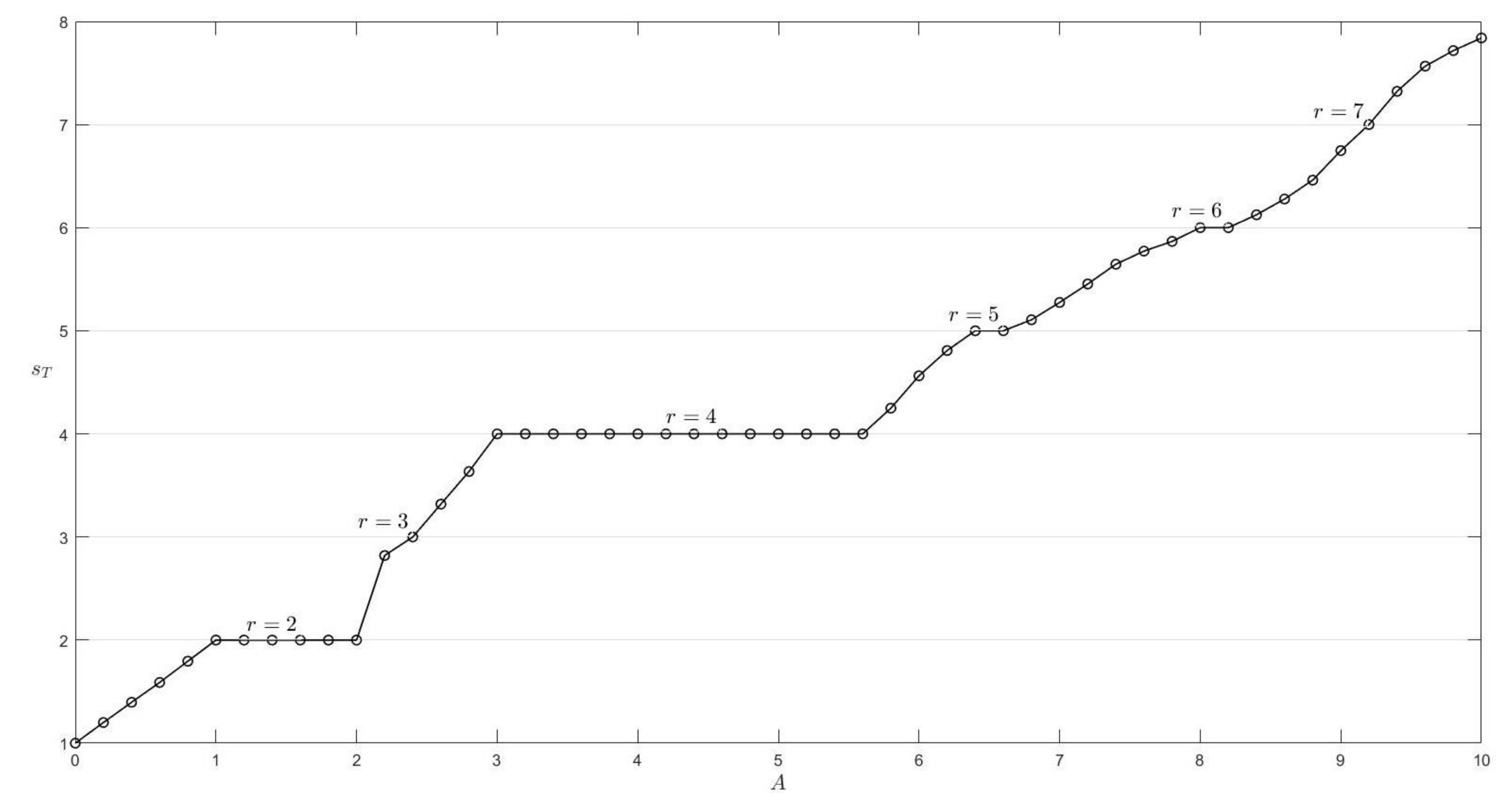}
\caption{Turbulent flame speed $s_T(A)$ for inviscid G-equation (\ref{Gi}) and unsteady cellular flow (\ref{Vu}) with $\omega=2$.}
\label{Fig_StGiVuA}
\end{figure}

\begin{figure}\center
\includegraphics[width=0.9\textwidth]{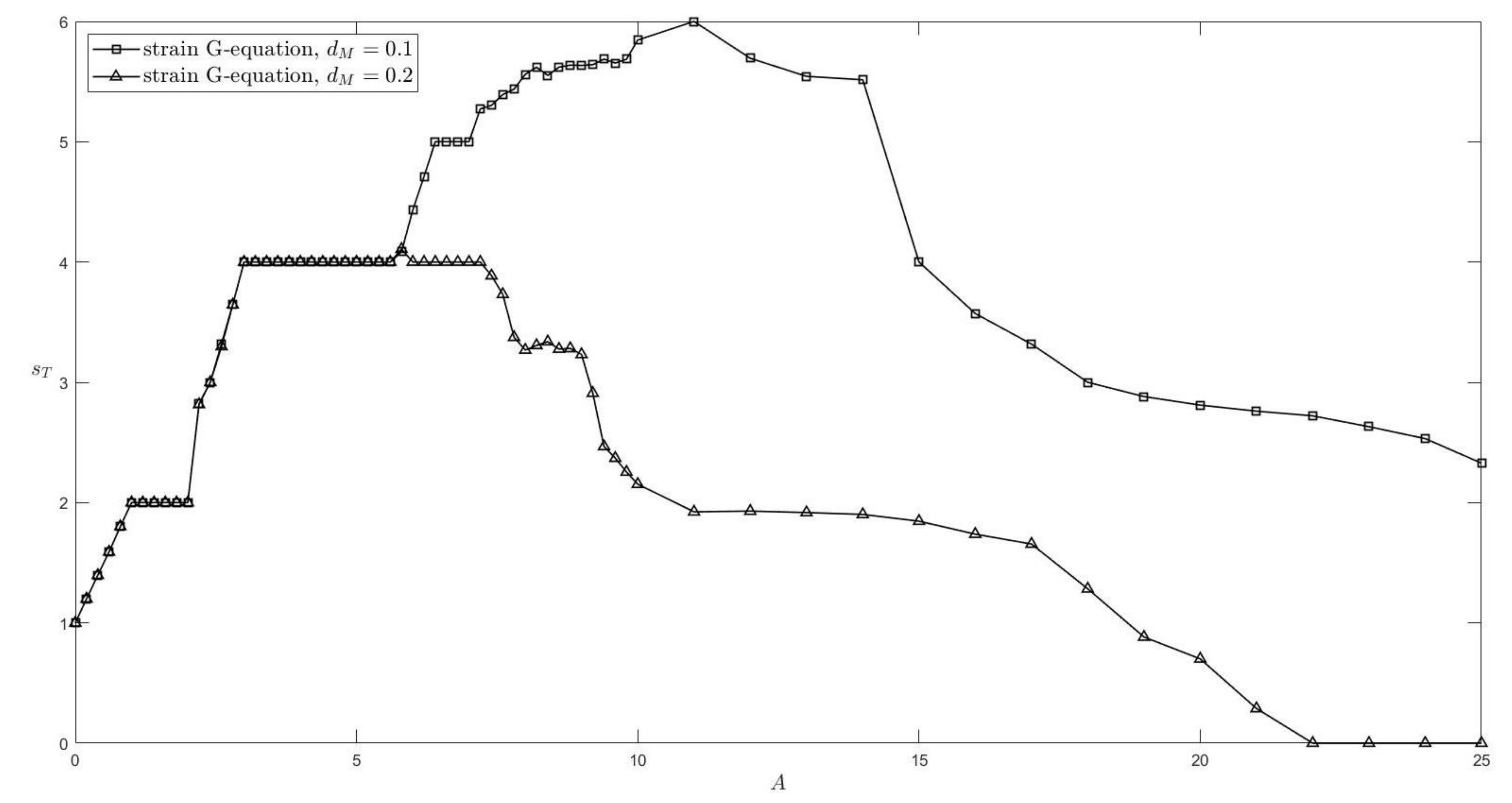}
\caption{Turbulent flame speed $s_T(A)$ for strain G-equation (\ref{Gs}) and unsteady cellular flow (\ref{Vu}) with $\omega=2$.}
\label{Fig_StGsVuA}
\end{figure}

Figure \ref{Fig_StGiVuA} and figure \ref{Fig_StGsVuA} are  plots of $s_T(A)$ for G-equations (\ref{Gi})(\ref{Gs}) and unsteady cellular flow (\ref{Vu}) with $\omega=2$. Two major differences are observed compare to figure \ref{Fig_StVsA} for steady cellular flow. As the front speed enhancement being synchronized with the time-periodic cellular flow, $s_T(A)$ becomes a piecewise constant function. Also as the flame front being disturbed by the unsteady cellular flow, flame quenching is delayed until much larger flow intensity.

\section{Conclusion}

We have performed a computational study on front propagation in G-equation models and the cellular flows. Two issues have been addressed in comparison to our previous works \cite{LXY13a,LXY13b}. First, a careful makeover of discretization of the strain rate is given so that the  monotonicity (with respect to one-sided derivatives) and hence the stability are met (even in absence of the curvature effect). Second, synchronization of front propagation may occur due to temporal oscillation in the Rayleigh-Bénard advection, and  the turbulent flame speeds may locally locked into a constant with respect to increase of the flow intensity.

Computation of turbulent flame speeds in G-equation models is rather challenging due to higher order discretization as well as large time simulation. The only exception so far might be the viscous G-equation as the curvature term being simplified to the diffusion term. In \cite{LXY13a}, turbulent flame speed is obtained as the effective Hamiltonian by solving the cell problem in homogenization theory. In \cite{GXZ21}, the viscous G-equation is discretized and simulated by the Galerkin proper orthogonal decomposition (POD) method. In future work, we plan to study accurate and efficient algorithms in solving G-equation models and evaluating turbulent flame speeds in three space dimensions.

Synchronization is a well-known nonlinear phenomena in chaotic dynamical systems. Besides having appeared in mathematical models like nonlinear oscillators or circle maps, synchronization has been widely applied in engineering science disciplines (for example, phase locking in circuit design). See \cite{AANVS06} for more details. In case of inviscid G-equation, its solutions are obtained by optimal control theory: $G(\x,t)=\inf_{\y(\cdot)} G(\y(t),0)$ with the infimum taken among all trajectories $\dot{\y}(\cdot)=\V(\y(\cdot),\cdot)+\a(\cdot)$, $\y(0)=\x$ and controls satisfying $|\a(\cdot)|\leq s_L$. Even the flow velocity has chaotic streamlines, the control effect given by the laminar velocity may contribute to self-organization of the trajectories. In future work, we would like to investigate the mechanism therein.

\medskip

\noindent {\bf Acknowledgments.} Yu-Yu Liu was supported by Ministry of Science and Technology grant 107-2115-M-006-017- of Taiwan. Jack Xin was partially supported by National Science Foundation grants DMS-1854434, DMS-1924548, DMS-1952644 of USA.

\begin{thebibliography}{99}

\bibitem{ACVV02}
M. Abel, M. Cencini, D. Vergni and A. Vulpiani.
{\em Front Speed Enhancement in Cellular Flows.}
Chaos: An Interdisciplinary Journal of Nonlinear Science 12 (2002), no. 2, 481--488.

\bibitem{AANVS06}
V. S. Anishchenko, V. Astakhov, A. Neiman, T. Vadivasova and L. Schimansky-Geier.
{\em Nonlinear Dynamics of Chaotic and Stochastic Systems.} Springer Berlin Heidelberg, 2006. 

\bibitem{BCVV95}
L. Biferale, A. Cristini, M. Vergassola and A. Vulpiani.
{\em Eddy Diffusivities in Scalar Transport.} 
Physics of Fluids 7 (1995), no. 11, 2725--2734.

\bibitem{B92}
D. Bradley.
{\em How Fast Can We Burn?} 
Symposium (International) on Combustion 24 (1992), no. 1, 247--262. 

\bibitem{CL84}
M. G. Crandall and P.-L. Lions.
{\em Two Approximations of Solutions of Hamilton-Jacobi Equations.} 
Mathematics of Computation 43 (1984), no. 167, 1--19. 

\bibitem{CNS11}
P. Cardaliaguet, J. Nolen and P. E. Souganidis.
{\em Homogenization and Enhancement for the G-Equation.}
Archive for Rational Mechanics and Analysis 199 (2011), no. 2, 527--561.

\bibitem{CTVV03}
M. Cencini, A. Torcini, D. Vergni and A. Vulpiani.
{\em Thin Front Propagation in Steady and Unsteady Cellular Flows.}
Physics of Fluids 15 (2003), no. 3, 679--688.

\bibitem{CW91}
R. Camassa and S. Wiggins.
{\em Chaotic Advection in a Rayleigh-Bénard Flow.} 
Physical Review A 43 (1991), no. 2, 774--797. 

\bibitem{GXZ21}
H. Gu, J. Xin and Z. Zhang.
{\em Error Estimates for a POD Method for Solving Viscous G-Equations in Incompressible Cellular Flows.} 
SIAM Journal on Scientific Computing 43 (2021), no. 1, A636--A662.

\bibitem{JP00}
G.-S. Jiang and D. Peng.
{\em Weighted ENO Schemes for Hamilton-Jacobi Equations.} SIAM Journal on Scientific Computing 21 (2000), no. 6, 2126--2143. 

\bibitem{LXY13a}
Y.-Y. Liu, J. Xin and Y. Yu.
{\em A Numerical Study of Turbulent Flame Speeds of Curvature and Strain G-equations in Cellular Flows.}
Physica D: Nonlinear Phenomena 243 (2013), no. 1, 20--31.

\bibitem{LXY13b}
Y.-Y. Liu, J. Xin and Y. Yu.
{\em Turbulent Flame Speeds of G-equation Models in Unsteady Cellular Flows.}
Mathematical Modelling of Natural Phenomena 8 (2013), no. 3, 198--205.

\bibitem{MK99}
A. J. Majda and P. R. Kramer.
{\em Simplified Models for Turbulent Diffusion: Theory, Numerical Modeling, and Physical Phenomena.}
Physics Reports 314 (1999), no. 4–5, 237--574. 

\bibitem{O00}
A. M. Oberman.
{\em Level Set Motion by Advection, Growth, and Mean Curvature as a Model for Combustion.}
Ph.D. Thesis, University of Chicago, USA, 2001.

\bibitem{OF02}
S. Osher and R. Fedkiw.
{\em Level Set Methods and Dynamic Implicit Surfaces.}
In Applied Mathematical Sciences, Springer New York, 2003. 

\bibitem{P00}
N. Peters.
{\em Turbulent Combustion.} 
Cambridge University Press, 2000.

\bibitem{R95}
P. D. Ronney.
{\em Some Open Issues in Premixed Turbulent Combustion.} 
In Modeling in Combustion Science (pp. 1--22), Springer Berlin Heidelberg, 1995. 

\bibitem{S07}
C.-W. Shu.
{\em High Order Numerical Methods for Time Dependent Hamilton-Jacobi Equations.}
In Lecture Notes Series, Institute for Mathematical Sciences, National University of Singapore (pp. 47--91), World Scientific, 2007. 

\bibitem{VCK03}
N. Vladimirova, P. Constantin, A. Kiselev, O. Ruchayskiy and L. Ryzhik.
{\em Flame Enhancement and Quenching in Fluid Flows.}
Combustion Theory and Modelling 7 (2003), no. 3, 487--508.

\bibitem{W85}
F.A. Williams.
{\em Turbulent Combustion.} 
In The Mathematics of Combustion (pp. 97--131), Society for Industrial and Applied Mathematics, 1985.

\bibitem{X00}
J. Xin.
{\em Front Propagation in Heterogeneous Media.} 
SIAM Review 42 (2000), no. 2, 161--230. 

\bibitem{X09}
J. Xin.
{\em An Introduction to Fronts in Random Media.} 
Springer New York, 2009. 

\bibitem{XY10}
J. Xin and Y. Yu.
{\em Periodic Homogenization of Inviscid G-equation for Incompressible Flows.}
Communications in Mathematical Sciences 8 (2010), no. 4, 1067--1078. 

\bibitem{XY13}
Y. Yu and J. Xin.
{\em Sharp Asymptotic Growth Laws of Turbulent Flame Speeds in Cellular Flows by Inviscid Hamilton-Jacobi Models.}
Annales de l’Institut Henri Poincaré C, Analyse Non Linéaire 30 (2013), no. 6, 1049--1068. 

\bibitem{XY14}
Y. Yu and J. Xin.
{\em Front Quenching in G-equation Model Induced by Straining of Cellular Flow.} 
Archive for Rational Mechanics and Analysis 214 (2014), no. 1, 1--34. 

\bibitem{ZCX15}
P. Zu, L. Chen and J. Xin,.
{\em A Computational Study of Residual KPP Front Speeds in Time-Periodic Cellular Flows in the Small Diffusion Limit.}
Physica D: Nonlinear Phenomena 311--312 (2015), 37--44. 

\end{thebibliography}
\end{document}